\documentclass[12pt,reqno]{amsart}

\usepackage{soul}
\usepackage{amsthm}
\usepackage{amssymb}
\usepackage{latexsym}
\usepackage{multicol}
\usepackage{verbatim,enumerate}
\usepackage{accents}
\usepackage[usenames]{color}
\usepackage{hyperref}
\usepackage{hyperref}
\usepackage{amsmath, amscd}
\usepackage{soul}
\usepackage{dsfont}
\usepackage{tikz,setspace}
\usetikzlibrary{matrix,arrows,decorations.pathmorphing}
\usetikzlibrary{positioning}

\advance\textwidth by 1.2in \advance\oddsidemargin by -.6in \advance\evensidemargin by -.6in

\theoremstyle{plain}
\newtheorem{prop}{Proposition}
\newtheorem{thm}{Theorem}
\newtheorem{lem}{Lemma}
\newtheorem{cor}{Corollary}

\theoremstyle{definition}

\theoremstyle{remark}

\newcommand{\lie}[1]{\mathfrak{#1}}

\newenvironment{pf}{\proof}{\endproof}
\newcounter{cnt}
 \makeatletter
\def\mydggeometry{\makeatletter\dg@YGRID=1\dg@XGRID=20\unitlength=0.003pt\makeatother}
\makeatother \theoremstyle{remark}

\numberwithin{equation}{section}
\makeatletter
\def\section{\def\@secnumfont{\mdseries}\@startsection{section}{1}%
  \z@{.7\linespacing\@plus\linespacing}{.5\linespacing}%
  {\normalfont\scshape\centering}}
\def\subsection{\def\@secnumfont{\bfseries}\@startsection{subsection}{2}%
  {\parindent}{.5\linespacing\@plus.7\linespacing}{-.5em}%
  {\normalfont\bfseries}}
\makeatother

\begin{document}

\title[A note on the fusion product decomposition of Demazure modules]{A note on the fusion product \\ decomposition of Demazure modules}
\author[R. Venkatesh and Sankaran Viswanath ]{R. Venkatesh and Sankaran Viswanath}\address{\noindent Department of Mathematics, Indian Institute of Science, Bangalore 5600112}
\email{rvenkat@iisc.ac.in}
\address{\noindent The Institute of Mathematical Sciences, HBNI, Chennai 600113}
\email{svis@imsc.res.in}
\thanks{RV is partially funded by the grants DST/INSPIRE/04/2016/000848, MTR/2017/000347, and an Infosys Young Investigator Award.
SV is partially funded by the grant MTR/2019/000071.}
\begin{abstract}
We settle the fusion product decomposition theorem for higher level affine Demazure modules for the cases
 $E^{(1)}_{6, 7, 8}, F^{(1)}_4$ and $E^{(2)}_{6}$, thus completing the main theorems of Chari et al. (J.
Algebra, 2016) and Kus et al. (Represent. Theory, 2016). We obtain a new combinatorial proof
for the key fact, that was used in Chari et al. (op cit.), to prove this decomposition theorem. 
We give a case free uniform proof for this key fact.
\end{abstract}

\maketitle


\section{Introduction}

Affine Demazure modules of higher level have been the subject of intensive study due to their connection with the representation theory of quantum affine algebras \cite{CV2015}, crystal bases \cite{Schilling},  the $X = M$ conjecture \cite{Naoi}, Macdonald polynomials \cite{ChariMacdonald} and more.

One of the key structural results concerning these modules is their decomposition into a fusion product of smaller Demazure modules of the same level. This {\em Steinberg type decomposition theorem}, proved in \cite{CSVW14, DV16, Venkat} (see also Theorem~\ref{mainthm} below) plays a vital role in understanding their graded structure and has many applications. For instance, this is used in \cite{Brito} to prove that certain level two affine Demazure modules in type $A$ coincide with the prime representations of the quantum affine algebra introduced by Hernandez-Leclerc \cite{HLMod} in the context of  monoidal categorification of cluster algebras. The decomposition theorem also serves as a base case for the study of fusion products of Demazure modules of unequal levels, and can be used to obtain defining relations for special modules of this kind \cite{Ravinder}.

However, the decomposition theorem was only proved for the untwisted affine algebras other than  $E^{(1)}_{6, 7, 8}$ and $F^{(1)}_4$ in \cite{CSVW14} and for twisted affine algebras other than $E^{(2)}_{6}$ in  \cite{DV16}. The obstruction lay in a technical result concerning the action of the affine Weyl group on weights, which was established in \cite[Proposition 3.5]{CSVW14} first for the untwisted affines of type $C$  and then by separate root sytem arguments in each of the types $A, B, D$ and $G$; additionally an appendix tabulated computational evidence in types $E$ and $F$.

In this note, we formulate a stronger version of this technical result (Proposition~\ref{mainresult}) and give a short, uniform proof for all affine types. We also point out in section~\ref{sec:rems} the new cases of the decomposition theorem and other key results of \cite{CSVW14} which now follow as consequences.

\section{Preliminaries}
\subsection{} We assume that the base field is complex numbers throughout the paper. 
  We refer to \cite{K90} for the general theory of affine Lie algebras.
We denote by $A$ an indecomposable affine Cartan matrix,  and by
$S$ the corresponding Dynkin diagram with the labeling of vertices as in Table Aff $1-3$ from \cite[pp. 54--55]{K90}. 
Let $\mathring{S}$ be the Dynkin diagram obtained from $S$ by dropping the  $0^{\mathrm{th}}$ node and let 
$\mathring{A}$ be the Cartan matrix, whose Dynkin diagram is $\mathring{S}$.
Let $\lie g$ and $\mathring{\lie g}$  be the  affine Lie algebra and the finite--dimensional simple Lie algebra associated to $A$ and $\mathring{A}$ , respectively. We shall realize $\mathring{\lie g}$ as a subalgebra of $\lie g$. 
We let $\lie h$ and $\mathring{\lie h}$ denote the Cartan subalgebras of $\lie g$ and $\mathring{\lie g}$ respectively, and let  $\lie h^{*}_\mathbb{R}$ and $\mathring{\lie h}^{*}_{\mathbb{R}}$ denote their real forms.

 We let $\Phi$ and $\mathring{\Phi}$ denote the sets of roots of $\lie g$ and $\mathring{\lie g}$  respectively.  We fix $\Pi=\{\alpha_0,\dots,\alpha_n\}$  a basis  for $\Phi$ such that $\mathring{\Pi}=\{\alpha_1,\dots,\alpha_n\}$ is a basis for $\mathring{\Phi}$. 
The weight lattice (resp. coweight lattice) of $\Phi$ is denoted by $P$ (resp. $P^\vee$) and the set of dominant integral weights is denoted by 
$P^+$. Similarly, the weight lattice (resp. coweight lattice) of $\mathring{\Phi}$ is denoted by $\mathring{P}$ (resp. $\mathring{P}^\vee$) and the set of dominant integral weights by $\mathring{P}^+$.

\subsection{Affine Weyl group}
We recall the key facts about affine Weyl groups following \cite{K90} and \cite{Carter}. Let $W$ and $\mathring{W}$ be the Weyl groups of  $\lie g$ and $\mathring{\lie g}$ respectively. We denote by $s_i$ the reflection associated to the simple root $\alpha_i$ for $0\le i \le n$.
Then we have $W=\langle s_0, s_1, \cdots, s_n \rangle$ and $\mathring{W}=\langle s_1, \cdots, s_n \rangle$.
Let $\langle\ , \rangle$ be the standard non-degenerate symmetric invariant form on $\lie g$. This determines a bijection $\lie h\to \lie h^{*}$. There is an action of $W$ on $\lie h$ and $\lie h^{*}$, compatible with this bijection.

For each $\alpha\in \mathring{\mathfrak{h}}^*$, we define $t_\alpha:\mathfrak{h}^*\to \mathfrak{h}^*$ by  $$t_{\alpha}(\Lambda)=\Lambda+\Lambda(c)\alpha - (\langle \Lambda, \alpha \rangle+\frac{\langle \alpha, \alpha \rangle}{2}\Lambda(c))\delta$$
for $\Lambda\in \lie h^{*}$. Here $c$ denotes the canonical central element of $\lie g$ and $\delta$ the unique indivisible positive imaginary root of $\Phi$.

Let $\theta$ denote the highest root of $\mathring{\Phi}$ and define
$a_0=2$ if  $\lie g$ is of type  $\tt A^{(2)}_{2n}\ (n\geq 1)$ and $a_0=1$ otherwise. Then 
$W=\mathring{W}\ltimes t_{M}$, where  $t_{M}=\{t_\mu: \mu\in M\}$ and $M$ is the sublattice of $\mathring{\lie h}^{*}_\mathbb{R}$ generated by the elements $w\left(\frac{1}{a_0}\theta\right)$ for all $w\in \mathring{W}.$ The explicit description of the lattice $M$ can be found in \cite[Page 414]{Carter}.
The affine Weyl group $W$ acts on the set $\lie h^{*}_{\mathbb{R}, 1}=\{\Lambda\in \lie h^{*}_\mathbb{R}: \Lambda(c)=1\}$ and this induces an action on the orbit space $\lie h^{*}_{\mathbb{R}, 1}/\mathbb{R}\delta \cong \mathring{\lie h}^{*}_\mathbb{R}$. The element $t_{\mu}\in M$ acts on $\mathring{\lie h}^{*}_\mathbb{R}$ as the translation $x \mapsto x+\mu$. Thus, $W$ acts as affine transformations on $\mathring{\lie h}^{*}_\mathbb{R}$.

Let $H_{\alpha^\vee, k} =\{\lambda \in \mathring{\lie h}^{*}:  \lambda(\alpha^\vee) = k\}$; this defines an affine hyperplane in $\mathring{\lie h}^{*}$ for $\alpha \in \mathring{\Phi}, k \in \mathbb{Z}$, where $\alpha^\vee$ is the coroot corresponding to the root $\alpha$.

Let $\mathring{\mathcal{H}}$ be the set of  hyperplanes $\left\{H_{\alpha^\vee, 0}: \alpha \in \mathring{\Phi} \right\}$.
The connected components of $\mathring{\lie h}^{*}_\mathbb{R}-\bigcup_{H\in \mathring{\mathcal{H}}}H$ are the Weyl chambers of $\mathring{\lie h}^{*}_\mathbb{R}$.
The elements of $\mathring{W}$ permute the hyperplanes in $\mathring{\mathcal{H}}$ and hence act on the set of Weyl chambers.
The set $C=\left\{\lambda\in\mathring{\lie h}^{*}_\mathbb{R}: \lambda(\alpha_i^\vee)> 0\ \text{for}\ 1\le i \le n\right\}$ is called the fundamental Weyl chamber.
The map $w\mapsto w(C)$ gives a bijection from  $\mathring{W}$ to the set of Weyl chambers and the closure $\bar{C}$ is a fundamental region for the $\mathring{W}$ action on $\mathring{\lie h}^{*}_\mathbb{R}$.

Let $\mathcal{H}$ be a set of affine hyperplanes $\left\{H_{\alpha^\vee, k}: \alpha \in \mathring{\Phi}, k\in Z_\alpha \right\} $ in $\mathring{\lie h}^{*}_\mathbb{R}$,  where
the sets $Z_\alpha$ are defined as follows (see \cite[Page 414]{Carter}):
$$Z_\alpha=\begin{cases}
2\mathbb{Z}& \text{if $\lie g$ is of type $ B^{(1)}_{n}\ (n\geq 2), C^{(1)}_{n}\ (n\geq 3), F^{(1)}_{4}$  and $\alpha$ is short}\\
3\mathbb{Z} & \text{if $\lie g$ is of type $G^{(1)}_{2}$ and $\alpha$ is short} \\
\frac{1}{2}\mathbb{Z}, & \text{if $\lie g$ is of type $A^{(2)}_{2n}\ (n\geq 1)$  and $\alpha$ is long }\\
\mathbb{Z}& \text{otherwise} \\
\end{cases}
$$
The connected components of $\mathring{\lie h}^{*}_\mathbb{R}-\bigcup_{H\in \mathcal{H}}H$ are called the alcoves of $\mathring{\lie h}^{*}_\mathbb{R}$. 
The elements of $W$ permute the  affine hyperplanes in $\mathcal{H}$ and hence act on the set $\mathcal{A}$ of alcoves.
The set 
$$A=\left\{\lambda\in\mathring{\lie h}^{*}_\mathbb{R}: \lambda(\alpha_i^\vee)> 0\ \text{for}\ 1\le i \le n, \lambda(\theta^\vee)<1/a_0\right\}$$ is called the fundamental alcove.  The map $w\mapsto w(A)$ gives a bijection from $W$ onto the set of alcoves $\mathcal{A}$. Moreover, the closure  $$\bar{A}=\left\{\lambda\in\mathring{\lie h}^{*}_\mathbb{R}: \lambda(\alpha_i^\vee)\ge 0\ \text{for}\ 1\le i \le n, \lambda(\theta^\vee)\le1/a_0\right\}$$ is a fundamental region for the $W$ action on $\mathring{\lie h}^{*}_\mathbb{R}$. 
Finally,  let $\Lambda_0 \in P^+$ be the fundamental weight corresponding to the $0^{\text{th}}$ vertex of the Dynkin diagram of $\lie g$ and  let $w_0$ denote the unique longest element in $\mathring{W}$.

\subsection{}
We will need the following well-known elementary lemma.
\begin{lem}\label{keylem}
Let $\mathcal{F}$ denote the set of alcoves contained in the fundamental Weyl chamber. Then  
  $\bar C=\bigcup_{B\in \mathcal{F}}\bar B$.
\end{lem}

\section{The technical result}\label{main}
Set $M^+=M\cap \bar{C}.$ The following is our main technical result which is crucial in proving Theorem \ref{mainthm}.
We refer to \cite[Section 4]{CSVW14} and \cite[Section 6]{DV16} for more details.
\begin{prop}\label{mainresult}
  Let $\lie g$ be an affine Lie algebra. 
Given $\ell\in \mathbb{N}$ and $\lambda\in \mathring{P}^+$, there exists $\mu\in M^+$ and $w\in \mathring{W}$ such that $wt_\mu(\ell\Lambda_0-\lambda)\in P^+$.
\end{prop}

\begin{onehalfspacing}
\begin{pf}
Let $\lambda'=-w_0(\lambda)$. Since $\lambda\in \mathring{P}^+$ and $-w_0(\mathring{P}^+)=\mathring{P}^+$, we have $\lambda' \in \mathring{P}^+$. Thus, the element ${\lambda'}/{\ell}\in \bar{C}.$ 
There exists $B\in \mathcal{F}$ such that ${\lambda'}/{\ell}\in \bar{B}$ by Lemma \ref{keylem}. 
Since $W$ acts simply transitively on the set of alcoves $ \mathcal{A}$, we have $\bar{B}=t_{\mu'} \,u\,\bar{A}$ for some ${\mu'}\in M$ and $u\in \mathring{W}$.
Since $0\in \bar{A}$, we have $\mu'=t_{\mu'} \,u(0)\in \bar B$. This implies that $\mu' \in \bar {B}\cap M\subseteq \bar{C}\cap M=M^+$.
Set $\mu=-w_0\, \mu'\in M^+$ and $w=u^{-1}w_0 \in \mathring{W}$ and
consider $$wt_{\mu}(\ell\Lambda_0-\lambda)\equiv \ell\Lambda_0+w(\ell\mu-\lambda)\ \text{mod}\ \mathbb{Z}\delta.$$ We claim that 
$wt_{\mu}(\ell\Lambda_0-\lambda)\in P^+.$ This will follow if we prove that (i) $w(\ell\mu-\lambda)\in  \mathring{P}^+$ and
(ii) $\ell\Lambda_0(\alpha_0^\vee)+w(\ell\mu-\lambda)(\alpha_0^\vee)\ge 0$.
Both these facts follow from the following observation: $$w(\ell\mu-\lambda)=u^{-1}(\lambda'-\ell\mu')=\ell u^{-1}\left(\frac{\lambda'}{\ell}-\mu'\right)=\ell u^{-1}t_{-\mu'}\left(\frac{\lambda'}{\ell}\right)\in \ell \bar{A}\cap \mathring{P}.$$ Since this belongs to $\ell \bar{A}$, we have  $w(\ell\mu-\lambda)(\theta^\vee) \leq \ell/a_0$. Recalling that $\alpha_0^\vee = c - a_0 \theta^\vee$ \cite[Chapter 6]{K90}, we conclude
$\ell + w(\ell\mu-\lambda)(\alpha_0^\vee) = \ell - a_0 w(\ell\mu-\lambda)(\theta^\vee)\ge 0$, proving (ii). Since $\ell \bar{A}\cap \mathring{P} \subseteq \mathring{P}^+$ we have
$w(\ell\mu-\lambda)\in \mathring{P}^+$, proving (i). 
\end{pf}
\end{onehalfspacing}

\section{Remarks}\label{sec:rems}
In the interest of completeness, we briefly point out the new cases of the main results of \cite{CSVW14} which hold, now that Proposition~\ref{mainresult} has been established.

\subsection{} As before, let $\mathring{\lie g}$ be a  finite-dimensional simple Lie algebra and $\lie g$ the corresponding untwisted affine Lie algebra.
We consider the Demazure modules of $\lie g$ that are stable under the current algebra $\mathring{\lie g}\otimes \mathbb{C}[t]$.
These modules are parametrized by pairs $(\ell, \lambda)$ with $\ell$ a positive integer and $\lambda$ a dominant integral weight of $\mathring{\lie g}$. Let $D(\ell, \lambda)$ denote the Demazure module corresponding to the pair $(\ell, \lambda)$.

\begin{thm}\label{mainthm}
  Let $\mathring{\lie g}$ be one of $E_6, E_7, E_8$ or $F_4$. Let $\lambda\in \mathring{P}^+$, $\ell \in \mathbb{N}$ and suppose that $\lambda=\ell \ (\sum_{i=1}^k \lambda_i)+\lambda_{0}$ with $\lambda_0 \in \mathring{P}^+$ and $\lambda_i \in {(\mathring{P}^\vee)}^+$ for $1 \leq i \leq k$. Then there is an isomorphism of $\mathring{\lie g}[t]$--modules
$$D(\ell,\lambda)\cong D(\ell,\ell\lambda_1)*D(\ell,\ell\lambda_2)*\cdots*D(\ell,\ell\lambda_k)* D(\ell,\lambda_0)$$
\end{thm}

\noindent
Here $*$ denotes the fusion product, and we refer to \cite[Section 2.7]{CSVW14} for the definition of fusion products of $\mathring{\lie g}\otimes \mathbb{C}[t]$--modules. As mentioned earlier, this was established in \cite[Theorem 1]{CSVW14} for all untwisted affines other than those of types $E, F$ and in \cite{DV16} for the twisted cases (modulo an extra hypothesis for $E_{6}^{(2)}$). It now also follows that \cite[Theorem 7]{DV16} holds without this extra condition.

\subsection{} Now let $\mathring{\lie g}$ be simply-laced, so that ${(\mathring{P}^\vee)}^+ = \mathring{P}^+$. Let $\alpha_i$ and $\omega_i$  for $1 \leq i \leq n$ denote the simple roots and fundamental weights respectively of $\mathring{\lie g}$. Given $\ell \geq 1$ and $\lambda \in \mathring{P}^+$, there is a unique decomposition $\lambda = \ell \, (\sum_{i=1}^n m_i \omega_i) + \lambda_0$ 
where $\lambda_0 = \sum_{i=1}^n r_i \omega_i$ with $0 \leq r_i <\ell$ and $m_i \geq 0$ for all $i$.
Theorem~\ref{mainthm} implies that for $\mathring{\lie g}$ of type $E$, one has
\begin{equation}\label{eq:primefact} 
  D(\ell,\lambda)\cong_{\mathring{\lie g}[t]}  D(\ell,\ell\omega_1)^{*m_1}*D(\ell,\ell\omega_2)^{*m_2}*\cdots*D(\ell,\ell\omega_n)^{*m_n}* D(\ell,\lambda_0)
\end{equation}

\medskip
Recall that a $\mathring{\lie g}[t]$-module is said to be {\em prime} if it is not isomorphic to a fusion product of non-trivial $\mathring{\lie g}[t]$ modules. 
In conjunction with \cite[Proposition 3.9]{CSVW14}, Theorem~\ref{mainthm} yields the following corollary:

\begin{cor}
Let $\mathring{\lie g}$ be one of $E_6, E_7$ or $E_8$. The decomposition \eqref{eq:primefact} gives a prime factorization of the Demazure module $D(\ell,\lambda)$, i.e., an expression as a fusion product of prime $\mathring{\lie g}[t]$-modules.
\end{cor}

This was established in \cite{CSVW14} for types $A$ and $D$.

\subsection{}
A second corollary concerns the notion of a $Q$-system \cite{Hatayama}. Roughly speaking a $Q$-system is a short exact sequence of $\mathring{\lie g}$-modules:
\[ 0 \rightarrow \bigotimes_{j\sim i} V(\ell \omega_j) \rightarrow V(\ell \omega_i) \otimes V(\ell \omega_i) \rightarrow  V((\ell+1) \omega_i) \otimes V((\ell-1) \omega_i) \rightarrow 0\]
where $\omega_i$ is a miniscule weight of $\mathring{\lie g}$ and $j\sim i$ means that nodes $i$ and $j$ are connected by an edge in the Dynkin diagram. Generalizations of $Q$-systems considered 
in \cite{CSVW14, CV2015, Fourier-Hernandez, DV16} involve replacing the tensor products above by fusion products of certain $\mathring{\lie g}[t]$-modules.
The following result  was established in \cite[Section 5]{CSVW14}:
\begin{prop}\label{prop:csvw-qsys} (\cite{CSVW14})
  Let $\mathring{\lie g}$ be simply-laced. Let $\ell \ge 1, \lambda \in \mathring{P}^+$ with $\ell \ge \max\{\lambda(\alpha^\vee): \alpha\in \mathring{\Phi}^+\}$ and suppose $\omega_i$ is a miniscule weight such that $\lambda(\alpha^\vee_i) >0$. Let $\mu = \ell \omega_i + \lambda - \lambda(\alpha^\vee_i) \alpha_i$. Then there exists a short exact sequence of $\mathring{\lie g}[t]$-modules:
\begin{equation}\label{eq:qsys-demzs}
0 \rightarrow \tau^*_{\lambda(\alpha^\vee_i)} \, D(\ell, \mu) \rightarrow D(\ell, \lambda + \ell \omega_i) \rightarrow D(\ell+1, \lambda + \ell \omega_i) \rightarrow 0
\end{equation}
\end{prop}

\medskip
For a graded module $V$, $\tau^*_dV$ denotes $V$ with its grading shifted by $d$; we refer to \cite{CSVW14} for a fuller explanation of the notations.
Proposition~\ref{prop:csvw-qsys} and Theorem~\ref{mainthm} together imply the following generalized $Q$-system in type $E$:
\begin{cor}
  Assume that $\mathring{\lie g}$ is of type $E_6$ or $E_7$, and retain the other notations of Proposition~\ref{prop:csvw-qsys}. Write $\mu=\ell \mu_1+\mu_0$ for some $\mu_1, \mu_0 \in \mathring{P}^+$. Then there exists a natural short exact sequence of $\mathring{\lie g}[t]$-modules:
\begin{gather*} 0\to\tau_{\lambda(\alpha^\vee_i)}^*\left( D(\ell, \ell\mu_1)*D(\ell, \mu_0)\right)\to D(\ell, \ell\omega_i)*D(\ell,\lambda)\\ \to D(\ell+1,(\ell+1)\omega_i)*D(\ell+1 ,\lambda-\omega_i)\to 0.\end{gather*}
\end{cor}

\medskip
We note that $\mathring{\lie g} = E_8$ is excluded since it does not have miniscule weights. This result was established in types $A, D$ in \cite[Theorem 2]{CSVW14}.




\end{document}